\documentclass[11pt, oneside, reqno]{amsart}
\usepackage{geometry}
\usepackage{graphicx}	
\usepackage{amssymb}

\newtheorem{prop}{Proposition}
\newtheorem{cor}{Corollary}

\newtheorem{conj}{Conjecture}

\title{Elemental Patterns from the Erd\H{o}s Straus Conjecture}
\author{Kyle Bradford}

\begin{document}
\maketitle
\begin{abstract}
This paper makes the following conjecture: For every prime $p$  there exists a positive integer $x$  with $\left\lceil \frac{p}{4} \right\rceil \leq x \leq \left\lceil \frac{p}{2} \right\rceil$  and a positive divisor $d|x^2$  so that either \\

\noindent (1) $ d \bmod \left( 4x - p \right) \equiv -px$ \\

\noindent or \\
 
\noindent (2) $d \leq x$  and $ d \bmod \left( 4x - p \right) \equiv -x$. \\

\noindent Furthermore this paper proves that the solutions to these modular equations are in one-to-one correspondence with the solutions of the diophantine equation used in the Erd\H{o}s Straus conjecture.
\end{abstract}

\section{Introductory Material}

 The Erd\H{o}s Straus conjecture suggests for any integer $n\geq2$ there exists positive integers $x,y$  and $z$  so that the following diophantine equation holds.
 
 \begin{equation} \label{eq: one}
 \frac{4}{n} = \frac{1}{x} + \frac{1}{y} + \frac{1}{z}
 \end{equation}

\ 

Introduced by Paul Erd\H{o}s and Ernst Straus in the late 1940s \cite{pe1}, the problem was quickly picked up by other notable mathematicians such as Richard Obl\'{a}th \cite{ob}, Luigi Rosati \cite{ro}, Koichi Yamamoto \cite{y} and Louis Mordell \cite{mo}.  Richard Guy included this problem in his book on Unsolved Problems in Number Theory along with many other results on Egyptian fractions \cite{rg}.  Notable papers use analytic number theory, abstractions or computational methods to analyze this problem \cite{aaa,bh,b,bi,ec,et,rg,iw,dl,gm,rav,san1,san2,san3,s,as,t,v,wweb1,wweb2,wweb3,y,yang}, but this paper introduces an insight that will govern how this problem will be resolved.  My earlier work described how it suffices to show the conjecture holds for any prime $p$ \cite{b}.  In this previous work I insisted that $x \leq y \leq z$ and I continue this convention. It was shown that $p \nmid x$,  $p|z$  and $p$  sometimes divides $y$.  It was also shown that $p^{2}$  does not divide $x, y$  or $z$.  Using the common nomenclature, a solution is of type I if $p \nmid y$  and is of type II if $p|y$.  This paper discusses new results and then motivates further research in the following sections.  All proofs are shown in the final section.

\section{New Results}

The results in this paper are very subtle but quite illuminating.  Ultimately, for each prime $p$, I build both necessary and sufficient conditions to describe the solutions of (\ref{eq: one})  solely through its smallest solution value, $x$.  The following proposition and corollary derive a sufficient condition for finding type I solutions to (\ref{eq: one}).

\begin{prop} \label{prop: one}
Suppose for a prime $p$  there exists a positive integer $x$  with $\left\lceil \frac{p}{4} \right\rceil \leq x \leq \left\lceil \frac{p}{2} \right\rceil$  and a positive divisor $d | x^2$  so that $ d \bmod \left( 4x - p \right) \equiv -px$. \\

\ 

\noindent Then letting

\begin{align*}
y &= \frac{px+d}{4x -p} \\
z &= \frac{p \left( x + p \left( \frac{x^{2}}{d} \right) \right)}{4x -p}
\end{align*}
\ 

\noindent we see that $x,y$  and $z$  are positive integers with $x \leq y \leq z$  and $p \nmid y$.
\end{prop}

\ 

\begin{cor} \label{cor: one}
Suppose for a prime $p$  there exists a positive integer $x$  with $\left\lceil \frac{p}{4} \right\rceil \leq x \leq \left\lceil \frac{p}{2} \right\rceil$  and a positive divisor $d|x^2$  so that $ d \bmod \left( 4x - p \right) \equiv -px$. \\

\ 

\noindent Then we have met a sufficient condition to find a type I solution to (\ref{eq: one}).
\end{cor}

\ 

The following proposition and corollary derive a sufficient condition for finding type II solutions to (\ref{eq: one}).

\begin{prop} \label{prop: two}
Suppose for a prime $p$  there exists a positive integer $x$  with $\left\lceil \frac{p}{4} \right\rceil \leq x \leq \left\lceil \frac{p}{2} \right\rceil$  and a positive divisor $d|x^2$  so that $d \leq x$  and $ d \bmod \left( 4x - p \right) \equiv -x$. \\

\ 

\noindent Then letting

\begin{align*}
y &=\frac{p(x+d)}{4x-p} \\
z &=\frac{p \left( x + \frac{x^{2}}{d} \right)}{4x -p}
\end{align*}
\ 

\noindent we see that $x,y$  and $z$  are positive integers with $x \leq y \leq z$  and $p|y$.
\end{prop}

\ 

\begin{cor} \label{cor: two}
Suppose for a prime $p$  there exists a positive integer $x$  with $\left\lceil \frac{p}{4} \right\rceil \leq x \leq \left\lceil \frac{p}{2} \right\rceil$  and a positive divisor $d|x^2$  so that $d \leq x$  and $ d \bmod \left( 4x - p \right) \equiv -x$. \\

\ 

\noindent Then we have met a sufficient condition to find a type II solution to (\ref{eq: one}).
\end{cor}

\ 

The following two propositions derive the necessary conditions for finding type I and type II solutions to (\ref{eq: one}) respectively.

\begin{prop} \label{prop: three}
Suppose for a prime $p$  there exist positive integers $x \leq y \leq z$ that satisfy (\ref{eq: one})  and $p \nmid y$. \\

\noindent Then it is necessarily true that $\left\lceil \frac{p}{4} \right\rceil \leq x \leq \left\lceil \frac{p}{2} \right\rceil$  and a positive divisor $d|x^2$  exists so that \\

\begin{align*}
(1) & \quad d \bmod \left( 4x - p \right) \equiv -px \\
(2) & \quad y = \frac{px+d}{4x -p}  \\
(3) & \quad z = \frac{p \left( x + p \left( \frac{x^{2}}{d} \right) \right)}{4x -p}.
\end{align*}
\end{prop}

\ 

\begin{prop} \label{prop: four}
Suppose for a prime $p$  there exist positive integers $x \leq y \leq z$ that satisfy (\ref{eq: one})  and $p|y$. \\

\noindent Then it is necessarily true that $\left\lceil \frac{p}{4} \right\rceil \leq x \leq \left\lceil \frac{p}{2} \right\rceil$  and a positive divisor $d|x^2$  exists so that \\

\begin{align*}
(1) & \quad d \leq x \\
(2) & \quad d \bmod \left( 4x - p \right) \equiv -x \\
(3) & \quad y = \frac{p(x+d)}{4x -p}  \\
(4) & \quad z = \frac{p \left( x +  \frac{x^{2}}{d} \right)}{4x -p}.
\end{align*}
\end{prop}

\ 

Indeed I have now developed both the necessary and sufficient conditions for solving the Erd\H{o}s Straus conjecture.  The solutions are in one-to-one correspondence with the modular identities in propositions \ref{prop: one}  and \ref{prop: two}.  This is a key result because it reduces the conjecture to one dimension.  That is to say, for every prime $p$  we need to find at least one pair, $x$  and $d$, meeting the appropriate conditions as functions of $p$ to prove the Erd\H{o}s Straus conjecture.  I summarize this in the following conjecture.  

\begin{conj} \label{conj:  one}
For every prime $p$  there exists a positive integer $x$  with $\left\lceil \frac{p}{4} \right\rceil \leq x \leq \left\lceil \frac{p}{2} \right\rceil$  and a positive divisor $d|x^2$  so that either \\

\noindent (1) $ d \bmod \left( 4x - p \right) \equiv -px$ \\

\noindent or \\
 
\noindent (2) $d \leq x$  and $ d \bmod \left( 4x - p \right) \equiv -x$.
\end{conj}

\section{Computation motivation toward a solution}

\begin{figure}[h] \label{fig: one}
\includegraphics{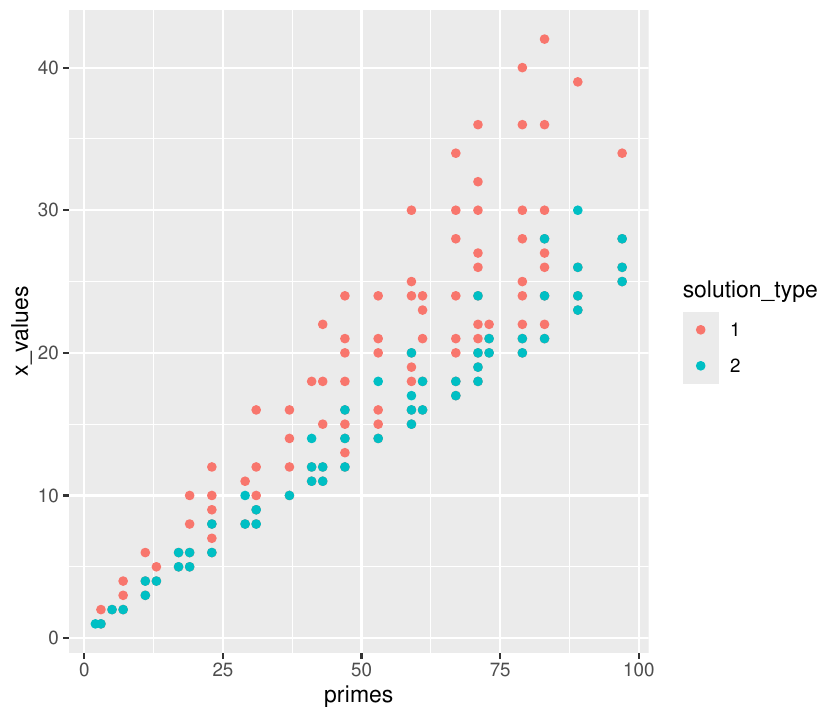}
\caption{This graphs primes less than 100 against the possible solution values x in the Erd\H{o}s Straus conjecture}
\end{figure}

The strength of this approach is that I have yet to employ different methodologies for different modular classes of prime numbers.  At this point you can use my conjecture to derive Mordell's identities for all primes $p$  except possibly for primes $p$  such that $p \bmod 840 \in \{ 1, 121, 169, 289,361,529 \}$.  It has been suggested that this problem can be solved through quadratic residues, and it may be no coincidence that my results suggest $d|x^{2}$.  My approach begs for $x$  to depend on a divisor of $\left\lceil \frac{p}{4} \right\rceil$, but this is clearly not the case for all primes.  Figure \ref{fig: one} may motivate you to find a similar pattern, although I only graphed primes less than 100 for clarity.

To make this more clear I include the following two tables.  In table \ref{tab: one} I consider type I solutions for $p < 100$  and I consider writing $x = \left\lceil \frac{p}{4} \right\rceil +k$.  Notice now that $0 \leq k \leq \left\lceil \frac{p}{4} \right\rceil$.   The table provides for each prime $p$ every possible value $k$  that appears for type I solutions.  Notice that for prime $p \neq 2$  so that $1 \bmod 24 \not\equiv p \bmod 24$  we are guaranteed to have a type I solution when $k=0$.  Also notice for primes $p$  so that $3 \bmod 4 \equiv p$  we are guaranteed to have type I solutions when $k$  is a divisor of $\left\lceil \frac{p}{4} \right\rceil$. \\

\ 
\begin{table}[h] \label{tab: one}
\begin{tabular}{|c|ccccccccccc|}
\hline
2 & & & & & & & & & & & \\
3 & 0 & 1 & & & & & & & & & \\
5 & 0 & & & & & & & & & & \\
7 & 0 & 1 & 2 & & & & & & & & \\
11 & 0 & 1 & 3 & & & & & & & & \\
13 & 0 & 1 & & & & & & & & & \\
17 & 0 & 1 & & & & & & & & & \\
19 & 0 & 1 & 3 & 5 & & & & & & & \\
23 & 0 & 1 & 2 & 3 & 4 & 6 & & & & & \\
29 & 0 & 3 & & & & & & & & & \\
31 & 0 & 1 & 2 & 4 & 8 & & & & & & \\
37 & 0 & 2 & 4 & 6 & & & & & & & \\
41 & 0 & 1 & 7 & & & & & & & & \\
43 & 0 & 1 & 4 & 7 & 11 & & & & & & \\
47 & 0 & 1 & 2 & 3 & 4 & 5 & 9 & 12 & & & \\
53 & 0 & 1 & 2 & 6 & 7 & 10 & & & & & \\
59 & 0 & 1 & 3 & 4 & 5 & 9 & 10 & 15 & & & \\
61 & 0 & 2 & 5 & 7 & 8 & & & & & & \\
67 & 0 & 1 & 3 & 4 & 7 & 11 & 13 & 17 & & & \\
71 & 0 & 1 & 2 & 3 & 4 & 6 & 8 & 9 & 12 & 14 & 18 \\
73 & 1 & 2 & 3 & & & & & & & \\
79 & 0 & 1 & 2 & 4 & 5 & 8 & 10 & 16 & 20 & & \\
83 & 0 & 1 & 3 & 5 & 6 & 7 & 9 & 15 & 21 & & \\
89 & 0 & 1 & 3 & 16 & & & & & & & \\
97 & 0 & 1 & 3 & 9 & & & & & & & \\
\hline
\end{tabular}
\ 

\ 

\caption{This table provides for each prime less than 100 all of the possible k values for type I solutions.}
\end{table}
\ 

In table \ref{tab: two} I consider type II solutions for primes $p<100$.  I again consider writing $x = \left\lceil \frac{p}{4} \right\rceil +k$  so that $0 \leq k \leq \left\lceil \frac{p}{4} \right\rceil$.  For a given prime $p$, notice that there are some type II solutions that have $x$  values that have no type I solutions.  For example, see that $p=41$  has $k=3$, which corresponds to $x=14$.  There are no type I solutions when $p=41$  and $x=14$. \\

\ 

\begin{table}[h] \label{tab: two}
\begin{tabular}{|c|cccc|}
\hline
2 & 0 & & & \\
3 & 0 & & & \\
5 & 0 & & & \\
7 & 0 & & & \\
11 & 0 & 1 & & \\
13 & 0 & & & \\
17 & 0 & 1 & & \\
19 & 0 & 1 & & \\
23 & 0 & 2 & & \\
29 & 0 & 2 & & \\
31 & 0 & 1 & & \\
37 & 0 & & & \\
41 & 0 & 1 & 3 & \\
43 & 0 & 1 & & \\
47 & 0 & 2 & 4 & \\
53 & 0 & 4 & & \\
59 & 0 & 1 & 2 & 5 \\
61 & 0 & 2 & & \\
67 & 0 & 1 & & \\ 
71 & 0 & 1 & 2 & 6 \\
73 & 1 & 2 & &  \\
79 & 0 & 1 & & \\
83 & 0 & 3 & 7 & \\
89 & 0 & 1 & 3 & 7 \\
97 & 0 & 1 & 3 & \\
\hline
\end{tabular}
\ 

\ 

\caption{This table provides for each prime less than 100 all of the possible k values for type II solutions.}
\end{table}
\ 

There are other patterns to consider.  My colleagues and I have considered tables like these for up to five digit primes, but this paper should outline a motivation for finding a pattern and proving the conjecture.

\section{Proofs}

\begin{proof}
\underline{Proposition \ref{prop: one}} \\

\noindent Let $p$  be a prime,  $x$  be a positive integer with $\left\lceil \frac{p}{4} \right\rceil \leq x \leq \left\lceil \frac{p}{2} \right\rceil$  and $d$  be a positive divisor $d|x^2$  so that $ d \bmod \left( 4x - p \right) \equiv -px$. \\

\noindent It should be clear that $p \neq 2$  because if $p=2$, then both $x$  and $d$  must be $1$  by definition and we cannot have $1 \bmod 2 \equiv 0$. \\

\noindent First, note that $x$  is a positive integer by definition. \\

\noindent  Next, note that if $ d \bmod \left( 4x - p \right) \equiv -px$, then $ (px+ d) \bmod \left( 4x - p \right) \equiv 0$.  This implies that $(4x-p)|(px+d)$.  By definition, $p,x$  and $d$  are all positive, so $px+d$  is positive.  Because $\left\lceil \frac{p}{4} \right\rceil \leq x \leq \left\lceil \frac{p}{2} \right\rceil$, it should be clear that $4x-p$  is also positive.  Letting $$y = \frac{px+d}{4x-p}$$
\ 

\noindent we see that $y$  is a positive integer.  We also see that $p \nmid y$.  To show this we assume, for the sake of contradiction, that $p|y$.  This implies that $p|(px+d)$  which further implies that $p|d$.  If $p|d$, then $p|x^{2}$;  however, from the definition of $x$,  we see that $x < p$.  $p$  cannot divide $x$  or $x^{2}$.  This creates a contradiction, implying that $p \nmid y$. \\

\noindent Finally, note that $\left( \frac{p^{2}x^{2}}{d} \right) d = p^{2}x^{2} = (-px)^{2}$  in $\mathbb{Z}$, so we get the following modular equation:
\begin{equation}
\left( \frac{p^{2}x^{2}}{d} \right) \bmod (4x-p) \cdot d \bmod (4x-p) \equiv (-px) \bmod (4x-p) \cdot (-px) \bmod (4x-p)
\end{equation}

\ 

\noindent Recall that $d \bmod (4x-p) \equiv -px$, so this equation becomes:

\begin{equation} \label{eq: two}
\left( \frac{p^{2}x^{2}}{d} \right) \bmod (4x-p) \cdot (-px) \bmod (4x-p) \equiv (-px) \bmod (4x-p) \cdot (-px) \bmod (4x-p)
\end{equation}

\ 

\noindent When $p \neq 2$, we see that $-px$  and $4x-p$  are coprime.  To show this let $m$  be a positive integer so that $m|(-px)$  and $m|(4x-p)$.  We necessarily see that $m$  must divide $(-4)(-px)+(-p)(4x-p) = p^{2}$.  This makes $m=1,p$  or $p^{2}$.  For sake of contradiction, assume that $m=p$.  This implies that $p|(4x-p)$,  which further implies that $p|x$, but this is impossible as $x \leq \left\lceil \frac{p}{2} \right\rceil$.  Next, for the sake of contradiction, assume that $m=p^{2}$.  This implies that $p^{2}|(-px)$, which further implies that $p|x$.  Again, this is impossible.  We conclude that $m=1$, so we have that $-px$  and $4x-p$  are indeed coprime.  We see that $(-px) \bmod (4x-p)$  is a unit, with an inverse element in the group $\left(\mathbb{Z} \slash (4x-p) \mathbb{Z} \right)^{\times}$.  Applying this inverse element on the right to either side of  (\ref{eq: two}), we have that $\left( \frac{p^{2}x^{2}}{d} \right) \bmod (4x-p) \equiv -px$.  We see then that $\left( px + \left( \frac{p^{2}x^{2}}{d} \right)   \right) \bmod (4x - p) \equiv 0$.  This implies that $(4x-p)|\left( px + \left( \frac{p^{2}x^{2}}{d} \right)   \right)$.  Note that $p,x,d$  and $4x-p$  are positive integers.  Letting $$ z =\frac{p \left( x + p \left( \frac{x^{2}}{d} \right) \right)}{4x -p}$$
\ 

\noindent we see that $z$  is a positive integer. \\

\noindent  To finish this proof we show that $x \leq y \leq z$.  \\

\noindent First consider when $ x < \left\lceil \frac{p}{2} \right\rceil$.  We see that $$x \leq \frac{p}{2} + \frac{d}{4x}$$

\ 

\noindent This implies that $4x^{2} < 2px +d$, $x(4x-p) < px+d$  and $$x < \frac{px+d}{4x-p} \leq y$$

\ 

\noindent Next consider when $ x = \left\lceil \frac{p}{2} \right\rceil$.  Because $p \neq 2$, we have that  $x = \frac{p+1}{2}$.  We see that 
\begin{align*}
y &= \frac{p \left( \frac{p+1}{2} \right) + d}{4 \left( \frac{p+1}{2} \right) - p} \\
&= \frac{p(p+1) + 2d}{2(p+2)} \\
&= \frac{p+1}{2} + \frac{d-(p+1)}{p+2}
\end{align*}

\ 

\noindent Because $y$  is an integer, we see that $(p+2)|(d - (p+1))$.  For $0<d<(p+1)$, we see that $$ -1 < \frac{d-(p+1)}{p+2} < 0$$

\ 

\noindent This implies that $d \geq (p+1)$  and $x = \frac{p+1}{2} \leq y$. \\

\noindent In either scenario we are guaranteed to have $x \leq y$. \\ 

\noindent  Because $d|x^{2}$, we see that $d \leq x^{2} \leq px$.  This implies that $d^{2} \leq p^{2}x^{2}$.  We see then that $$ d \leq \left( \frac{p^{2}x^{2}}{d} \right)$$

\ 

\noindent We see then that $$ \frac{px + d}{4x- p} \leq \frac{px + \left( \frac{p^{2}x^{2}}{d} \right)}{4x-p} $$

\ 

\noindent We see that $y \leq z$.
\end{proof}

\ 

\begin{proof}
\underline{Corollary \ref{cor: one}} \\

\noindent Let $p$  be prime. \\

\noindent Suppose that a positive integer $x$  exists with $\left\lceil \frac{p}{4} \right\rceil \leq x \leq \left\lceil \frac{p}{2} \right\rceil$  and a positive divisor $d|x^2$  exists so that $ d \bmod \left( 4x - p \right) \equiv -px$. \\

\noindent Under these conditions, we see from Proposition \ref{prop: one} that positive integers $x,y$  and $z$  exist with $x \leq y \leq z$  and $p \nmid y$  so that

\begin{align*}
\frac{1}{x} + \frac{1}{y} + \frac{1}{z} &= \frac{1}{x} + \frac{4x-p}{px+d } + \frac{4x-p}{p \left(x+p \left( \frac{x^{2}}{d} \right) \right)} \\
&= \frac{p \left( px + d \right) \left( x + p \left(  \frac{x^2}{d} \right) \right) + px\left(4x-p \right) \left(x + p\left( \frac{x^2}{d} \right) \right) + x \left( px+d \right) \left( 4x - p \right)}{px \left( px + d \right) \left(x + p \left( \frac{x^2}{d} \right) \right)} \\
&= \frac{2p^2x^2 + \frac{p^3x^3}{d} + dpx + 4px^3 + \frac{4p^2x^4}{d} - p^2x^2 - \frac{p^3x^3}{d}+4px^3+4dx^2 - p^2x^2 -dpx}{p \left( 2px^3 + \frac{p^2 x^4}{d} +dx^2 \right)} \\
&= \frac{4 \left(2px^3 + \frac{p^2x^4}{d}+dx^2 \right)}{p \left( 2px^3 + \frac{p^2x^4}{d}+dx^2  \right)} \\
&= \frac{4}{p}
\end{align*}

\ 

\noindent This is a type I solution to (\ref{eq: one}).
\end{proof}

\ 

\begin{proof}
\underline{Proposition \ref{prop: two}} \\

\noindent Let $p$  be a prime,  $x$  be a positive integer with $\left\lceil \frac{p}{4} \right\rceil \leq x \leq \left\lceil \frac{p}{2} \right\rceil$  and $d$  be a positive divisor $d | x^2$  so that $d \leq x$  and $ d \bmod \left( 4x - p \right) \equiv -x$. \\

\noindent First, note that $x$  is a positive integer by definition. \\

\noindent  Next, note that if $ d \bmod \left( 4x - p \right) \equiv -x$, then $ (x+ d) \bmod \left( 4x - p \right) \equiv 0$.  This implies that $(4x-p)|(x+d)$.  By definition, $p,x,d$  and $4x-p$  are all positive, so $p(x+d)$  is positive.  Letting $$y = \frac{p(x+d)}{4x-p}$$
\ 

\noindent we see that $y$  is a positive integer.  Because $(4x-p)|(x+d)$,  we see that $p|y$  in that scenario.  \\

\noindent Finally, note that $\left( \frac{x^{2}}{d} \right) d = x^{2} = (-x)^{2}$  in $\mathbb{Z}$, so we get the following modular equation: 
\begin{equation}
\left( \frac{x^{2}}{d} \right) \bmod (4x-p) \cdot d \bmod (4x-p) \equiv (-x) \bmod (4x-p) \cdot (-x) \bmod (4x-p)
\end{equation}
\ 

\noindent Recall that $d \bmod (4x-p) \equiv -x$, so this equation becomes:
\begin{equation} \label{eq: four}
\left( \frac{x^{2}}{d} \right) \bmod (4x-p) \cdot (-x) \bmod (4x-p) \equiv (-x) \bmod (4x-p) \cdot (-x) \bmod (4x-p)
\end{equation}
\ 

\noindent  We see that $-x$  and $4x-p$  are coprime.  To show this let $m$  be a positive integer so that $m|(-x)$  and $m|(4x-p)$.  We necessarily see that $m$  must divide $(-4)(-x)+(-1)(4x-p)=p$.  This makes $m=1$  or $p$.  For sake of contradiction assume that $m=p$.  This implies that $p|(-x)$, which further implies that $p|x$.  We have shown this to be impossible.  We conclude that $m=1$, so we have that $-x$  and $4x-p$  are indeed coprime.  We see that $ (-x) \bmod (4x-p)$  is a unit, with an inverse element in the group $\left(\mathbb{Z} \slash (4x-p) \mathbb{Z} \right)^{\times}$.  Applying this inverse element on the right to either side of (\ref{eq: four}), we have that $\left( \frac{x^{2}}{d} \right) \bmod (4x-p) \equiv -x$.  We see then that $\left( x + \left( \frac{x^{2}}{d} \right)   \right) \bmod (4x - p) \equiv 0$.  This implies that $(4x-p)|\left( x + \left( \frac{x^{2}}{d} \right)   \right)$.  Note that $p,x,d$  and $4x-p$  are positive integers.  Letting $$ z =\frac{p \left( x + \left( \frac{x^{2}}{d} \right) \right)}{4x -p}$$
\ 

\noindent we see that $z$  is a positive integer.  \\

\noindent To finish this proof we show that $x \leq y \leq z$. \\

\noindent  Because $p|y$, we see that $p \leq y$.  By definition we see that $x \leq p$.  This implies that $x \leq y$.  \\

\noindent  Because $d \leq x$, we have $d^{2} \leq x^{2}$.  We see that $$ d \leq \left( \frac{x^{2}}{d} \right)$$

\ 

\noindent which implies that $$ \frac{p(x + d)}{4x- p} \leq \frac{p \left(x + \left( \frac{x^{2}}{d} \right) \right)}{4x-p} $$

\ 

\noindent We see that $y \leq z$.  
\end{proof}

\ 

\begin{proof}
\underline{Corollary \ref{cor: two}} \\

\noindent Let $p$  be prime.  \\

\noindent Suppose that positive integer $x$  with $\left\lceil \frac{p}{4} \right\rceil \leq x \leq \left\lceil \frac{p}{2} \right\rceil$  and a positive divisor $d|x^2$  so that $d \leq x$  and $ d \bmod \left( 4x - p \right) \equiv -x$. \\

\noindent Under these conditions, we see that positive integers $x,y$  and $z$  exist with $x \leq y \leq z$  and $p | y$  so that

\begin{align*}
\frac{1}{x} + \frac{1}{y} + \frac{1}{z} &= \frac{1}{x} + \frac{4x-p}{p(x+d)} + \frac{4x-p}{p \left( x + \frac{x^{2}}{d} \right)} \\
&= \frac{p(x+d) \left( x + \frac{x^{2}}{d} \right) + x(4x-p)\left( x + \frac{x^{2}}{d} \right) + x(x+d)(4x-p)}{px(x+d) \left( x + \frac{x^{2}}{d}\right)} \\
&= \frac{2px^2 + \frac{px^3}{d} +dpx + 4x^3 + \frac{4x^4}{d} -px^2 - \frac{px^3}{d} +4x^3 + 4dx^2 - px^2 - dpx}{p \left( 2x^3 + \frac{x^4}{d} + dx^2 \right)} \\
&= \frac{4\left( 2x^3 + \frac{x^4}{d} + dx^2 \right)}{p \left( 2x^3 + \frac{x^4}{d} + dx^2 \right)} \\
&= \frac{4}{p}
\end{align*}

\end{proof}

\ 

\begin{proof}
\underline{Proposition \ref{prop: three}} \\

\noindent Let $p$  be prime and let $x \leq y \leq z$  be positive integers that satisfy (\ref{eq: one}) with $p \nmid y$.  \\

\noindent First, it was shown in \cite{b}  that a necessary condition for type I solutions is that $\left\lceil \frac{p}{4} \right\rceil \leq x \leq \left\lceil \frac{p}{2} \right\rceil$.  \\

\noindent Next, slightly changing the notation in \cite{b}, we let $m= \gcd(x,y,z), a=\gcd(x,y)/m, b=\gcd(x,z)/m$  and $c=\gcd(y,z)/m$.  Using this new notation, it was shown in \cite{b}  that $x=abm$, $y=acm$  and $z=bcm$. \\

\noindent For this type I solution, let $d= (4x-p)y - px$.  We need to show that $d|x^{2}$. \\

\noindent In \cite{b}  it was shown for type I solutions that $$ p = \frac{4abcm - a}{b+c}.$$
\ 

\noindent  This makes 
\begin{align*}
d &= (4x-p)y - px \\
&= 4xy -(x+y)p \\
&= 4a^{2}bcm^{2} - (abm + acm)\left( \frac{4abcm - a}{b+c} \right) \\
&= 4a^{2}bcm^{2} - 4a^{2}bcm^{2} + a^{2}m \\
&=a^{2}m.
\end{align*}
\ 

\noindent Because $a^{2}m | a^{2}b^{2}m^{2}$, we see that $d|x^{2}$. \\

\noindent From the definition of $d = (4x-p)y - px$, it should be clear that $d \bmod (4x-p) \equiv -px$. \\

\noindent We see that
\begin{align*}
\frac{px+d}{4x-p} &= \frac{px + (4x-p)y - px}{4x-p} \\
&= \frac{(4x-p)y}{4x-p} \\
&= y.
\end{align*}
\ 

\noindent We also see that
\begin{align*}
\frac{p\left(x+ p \left( \frac{x^{2}}{d} \right) \right)}{4x-p} &= \frac{px + \frac{p^{2}x^{2}}{(4x-p)y - px}}{4x-p} \\
&= \frac{(4x-p)xyp-p^{2}x^{2} + p^{2}x^{2}}{(4x-p)(4xy-(x+y)p)} \\
&= \frac{(4x-p)xyp}{(4x-p)(4xy-(x+y)p)} \\
&= \frac{xyp}{4xy-(x+y)p} \\
&= \frac{1}{\frac{4}{p} - \frac{1}{x} - \frac{1}{y}} \\
&= z.
\end{align*}
\end{proof}

\ 

\begin{proof}
\underline{Proposition \ref{prop: four}} \\

\noindent Let $p$  be prime and let $x \leq y \leq z$  be positive integers that satisfy (\ref{eq: one}) with $p|y$.  \\

\noindent It was shown in \cite{b}  that a necessary condition for type II solutions is that $\left\lceil \frac{p}{4} \right\rceil \leq x \leq \left\lceil \frac{p}{2} \right\rceil$.  \\

\noindent For this type II solution, let $d= (4x-p)(y \slash p) - x$.  We need to show that $d|x^{2}$. \\

\noindent Recall, using the slightly different notation, that $p|y$  and $p|z$,  so $p|c$.  In \cite{b}  it was shown for type II solutions that if $c = c^{*}p$, then $$ p = 4abm - \frac{a+b}{c^{*}}.$$
\ 

\noindent  This makes 
\begin{align*}
d &= (4x-p)(y \slash p) - x \\
&= 4x(y \slash p) -(x+y) \\
&= 4a^{2}bc^{*}m^{2} - abm - ac^{*}m \left(  4abm - \frac{a+b}{c^{*}} \right) \\
&= 4a^{2}bc^{*}m^{2} - abm - 4a^{2}bc^{*}m^{2} + a^{2}m + abm \\
&=a^{2}m.
\end{align*}
\ 

\noindent Because $a^{2}m | a^{2}b^{2}m^{2}$, we see that $d|x^{2}$. \\

\noindent We have that $y \leq z$  implies that $acm \leq bcm$  or $a \leq b$.  This implies that $a^{2}m \leq abm$,  so $d \leq x$. \\

\noindent From the definition of $d = (4x-p)(y \slash p) - x$, it should be clear that $d \bmod (4x-p) \equiv -x$. \\

\noindent We see that
\begin{align*}
\frac{p(x+d)}{4x-p} &= \frac{p (x + (4x-p)(y\slash p) - x)}{4x-p} \\
&= \frac{(4x-p)y}{4x-p} \\
&= y.
\end{align*}
\ 

\noindent We also see that
\begin{align*}
\frac{p\left(x+ \frac{x^{2}}{d} \right)}{4x-p} &= \frac{p \left(x + \frac{x^{2}}{(4x-p)(y \slash p) - x} \right)}{4x-p} \\
&= \frac{(4x-p)xy-px^{2} + px^{2}}{(4x-p)(4x(y \slash p)-(x+y))} \\
&= \frac{(4x-p)xy}{(4x-p)(4x(y \slash p)-(x+y))} \\
&= \frac{xy}{4x(y \slash p)-(x+y)} \\
&= \frac{1}{\frac{4}{p} - \frac{1}{x} - \frac{1}{y}} \\
&= z.
\end{align*}
\end{proof}







\bibliographystyle{amsplain}

\end{document}